\documentclass[a4paper]{amsart}
\usepackage{amssymb}
\usepackage{verbatim}
\usepackage{amscd}

\numberwithin{equation}{section}
\theoremstyle{plain}
\newtheorem{theorem}[equation]{Theorem}

\newtheorem{proposition}[equation]{Proposition}
\newtheorem{lemma}[equation]{Lemma}
\newtheorem*{(DQ1)}{(DQ1)}
\theoremstyle{definition}

\theoremstyle{remark}

\begin{document}
\title{On certain subshifts and their associated monoids}
\author{Toshihiro Hamachi}
\author{Wolfgang Krieger}
\begin{abstract} Within a subclass of monoids (with zero) a structural characterization is given of those that are 
associated to topologically transitive subshifts with Property (A). 

 \end{abstract}
 
\maketitle

\section{Introduction}
Let $\Sigma$ be a finite alphabet, and let $S_{\Sigma}$ be the shift 
on the shift space $\Sigma^{\Bbb Z}$,
$$
S_{\Sigma}((x_{i})_{i \in \Bbb Z}) =  (x_{i+1})_{i \in \Bbb Z}, 
\qquad 
(x_{i})_{i
\in \Bbb Z}  \in \Sigma^{\Bbb Z}.
$$
In symbolic dynamics one studies subshifts, by which are meant the dynamical systems $(C, S)$, where  
$C$ is an $S_{\Sigma}$-invariant closed subset of $\Sigma^{\Bbb Z}$, and $S$ is
the restriction of $S_{\Sigma}$ to $C$. An introduction to the theory 
of subshifts is in 
 \cite {Ki, LM} (see also \cite {B, BP}). 
 Here we continue the 
investigation of the semigroup invariants of subshifts as begun 
in 
 \cite {Kr}.
 There a partially ordered set $\mathcal P_{C}$  was invariantly 
associated to a subshift. Also an invariant property $(A)$ of 
subshifts was described. If a subshift has Property $(A)$  
then a semigroup (with zero) $\mathcal S_{C}$  is invariantly associated 
to $ C$. The semigroup  $\mathcal S_{C}$  is a monoid if and only if $\mathcal P_{C}$ is a
singleton set and then $\mathcal P_{C}$  contains the unit of $\mathcal S_{C}$. 

Let $\mathcal S$ be a semigroup (with zero). We denote by 
$\Gamma_{-}(a), \Gamma(a),\Gamma_{+}(a)$ the left 
context, context, and right context (in $\mathcal S$) of an element $a $ of $\mathcal 
S$,
\begin{align*}
&\Gamma_{-}( a )
 = \{ g^{-} \in \mathcal S :  g^{-}  a \neq 0 \},\\
 &
\Gamma( a ) = \{( g^{-}   ,  g ^{+}   ) \in \mathcal S   
\times  \mathcal S
 : g^{-}a g^{+}\neq 0 \},        \\
 &
\Gamma_{+}( a )
 = \{ g^{+} \in \mathcal S :    a g ^{+}\neq 0 \}. 
\end{align*}
Denoting by $[a]$ the set of elements in $\mathcal S$  that have the 
same context as $a\in \mathcal S$, the set  $[\mathcal S] = \{[a] : a  \in \mathcal 
S\}$  is  a semigroup (with zero), 
where
$$
[a][b] = [a b ], \qquad a , b  \in \mathcal S,
$$
(e.g. see \cite[Section 2.2]{P}).
We define subsemigroups $\mathcal S^-$  and  $\mathcal S^+$ of $\mathcal S$ by 
$$
\mathcal S^- = \bigcap_{b \in \mathcal S \setminus \{ 0 \}}
 \{a \in \mathcal S: ba   \neq 0    \},
\quad
\mathcal S^+ = \bigcap_{b \in \mathcal S \setminus \{ 0 \}}
 \{a \in \mathcal S:a   b \neq 0    \}.
$$

In section 2 we show  for a topologically transitive subshift $ C$ 
with property $(A)$, that the projection of $\mathcal  S_{C}$ onto $[\mathcal S_{C}]$ is an isomorphism. We also show  for a 
topologically transitive subshift $C$ with property $(A)$ such 
that 
$$
\mathcal  S ^{-}_{C} \cap \mathcal  S ^{+}_{C}   \neq \emptyset,
$$
that $\mathcal  S_{C}$ is a monoid (with zero), and that
$$
\mathcal  S ^{-}_{C} \cap  \mathcal  S ^{+}_{C}   = \{ \bold 1 \},
$$
or, equivalently, that
$$
0 \in( \bigcap_{a^{+}\in \mathcal S^+_{C}\setminus \{\bold 1 \}}\mathcal S^-_{C}a^{+} )
\cap( \bigcap_{a^{-}\in \mathcal S^-_{C}\setminus \{\bold 1 \}}a^{-}\mathcal S^+_{C} ) .    
$$

In the converse direction we consider in section 3 finitely generated 
monoids $\mathcal  M$ (with zero) such that
\begin{align*}
\mathcal  M \setminus \{0 \} =\mathcal  M ^{+}  \mathcal  M ^{-}  , \tag {$\star$}
\end{align*}
and such that
\begin{align*}
\{0, \bold 1\}  \subset( \bigcap_{a^{+}\in \mathcal M^+\setminus \{\bold 1  \}}\mathcal M^-
 a^{+})
\cap( \bigcap_{a^{-}\in \mathcal M^-\setminus \{\bold 1  \}}a^{-}\mathcal M^+ ) . \tag {$\star$$\star$}
\end{align*}
For a finite generating set  $\Delta 
$  
of a  monoid $\mathcal  M$ such that  ({$\star$})  and ({$\star$$\star$})  hold, we show that the topologically transitive subshift
$
C(\Delta) \subset \Delta^{\Bbb Z }
$
with admissible words $( \delta _{i})_{ 1 \leq i \leq I}\in \Delta^{[1,I]}, I \in \Bbb N,$   given by
$$
\prod _{ 1 \leq i \leq I}  \delta _{i} \ne 0,
$$
has  property $(A)$, and that 
$$
\mathcal  S_{C(\Delta)}= [\mathcal  M].
$$
As a consequence, a finitely generated monoid (with zero) $\mathcal  M$, such that ($\star$) 
and ($\star$$\star$) hold, is associated to a topologically transitive subshift with Property $(A)$ if and only if the projection of $\mathcal  M$ onto  $[\mathcal  M]$ is an isomorphism.

In section 4 we give examples. These we describe by confluent 
systems of string-rewriting rules.

\section{Topologically transitive subshifts with Property $(A)$}

 	We fix notation and terminology. Given a subshift $C \subset  
\Sigma^{ \Bbb Z}$, we set for $ x \in C$
$$
x_{[i, k]} = (x_{i})_{i \leq j \leq k },
$$
and we set
$$
C_{[i,k]} =\{x_{[i, k]} :x \in C \}, \qquad i, k \in \Bbb Z, i \leq k.
$$
We use similar notation also for blocks,
$$
b_{[i^{\prime}, k^{\prime}]} = (b_{j})_{i^{\prime} \leq j \leq k^{\prime} },
\qquad b \in C_{[i, k ]}, \  i\leq i^{\prime} \leq k^{\prime} \leq k , \ \ i, k \in
\Bbb Z.
$$
The left context of a block  $a\in C_{[i,k]}, i, k \in \Bbb Z, i \leq  k$, is 
denoted by $\Gamma^-(a)$,
$$
\Gamma ^{-}(a) = \bigcup_{J\in \Bbb N} \{ b\in C_{   [ i-J, i)}: ( b   , a      ) \in C_{  [ i-J, k]} \}.
$$
The notation $\Gamma^{+}(a)$ has the time symmetric meaning.
The context of a block 
 $a\in C_{[i,k]}, i, k \in \Bbb Z, i \leq  k$, is 
denoted by $\Gamma(a)$,
$$
\Gamma (a) =
\bigcup_{J, M \in \Bbb N} \{ (b, c)\in C_{   [ i-J, i)}\times C_{  (k,  k+M] }
: ( b   , a ,c     ) \in C_{  [ i-J, k+M] } \}.
$$
We also set
\begin{align*}
\omega^{+} (a) = \bigcap _{b\in \Gamma _{-}(a)}
\{ c \in\Gamma^+(a ):     b         \in  \Gamma^-(a,c)\},  \quad
 a \in
C_{[i,k]}, \  i,k \in \Bbb Z, \ i \leq k.
\end{align*}
The notation $\omega^{-}(a)$ has the time symmetric meaning. 

The  (smallest) period of a periodic point $p$ of $C$ we denote  by $\Pi(p)$.

We recall from [Kr] the construction of the partially ordered set $\mathcal  P_{C}$. For a 
subshift $C \subset  \Sigma^{ \Bbb Z}$ set
$$
X^{-}_{n} (C) = \bigcap _{ i \in \Bbb Z} \{ x \in C: x_{i} \in \omega ^{+}
(x_{[i-n, i)}) \}, \qquad n \in \Bbb N,
$$
with the time symmetric meaning for $ X^{+}_{n}(C)$. We set
$$
X_{n} (C) = X^{-}_{n} (C) \cap X^{+}_{n} (C),  \qquad
 n \in \Bbb N,
$$ Also set
$$
X^{+}(C) = \bigcup_{ n \in \Bbb N } X^{+}_{n} (C).
$$
The notation  $ X^{-}(C)$  has the time symmetric meaning. Also set
$$
X(C) =  \bigcup_{ n \in \Bbb N } X_{n} (C).
$$
In this paper we consider subshifts $C $ such that $X(C)$ is dense in $C$.
We denote by $P(X_{n}(C))$ the set of periodic points in $X_{n}(C), n\in 
\Bbb N$, and also set
$$
P(X(C) )= \bigcup _{n \in \Bbb N } P(X_{n}(C)).
$$
Introduce a reflexive and transitive relation $\lesssim_C$ into the set 
$P(X(C))$. For $u, u'\in P(X(C)), u \gtrsim_C u^{\prime}$ will mean that 
there exists a point in $X(C)$  that is negatively asymptotic to the 
orbit of $u$ and positively asymptotic to the orbit of $u'$. Denote 
the resulting equivalence relation by $\sim_C$ and the set of 
$\sim_C$-equivalence classes by  $\mathcal  P_{C}$. The partial 
order that is induced by $\lesssim_C$ on $\mathcal  P_{C}$ is also denoted by 
$ \lesssim_C$. ²

We recall from \cite {Kr} the definition of Property (A).
For $n \in \Bbb N$ a subshift $C \subset  \Sigma ^{\Bbb 
Z},$ 
has property $(a, n, H), H \in \Bbb N,$  if for $h,
\widetilde {h} \geq 3H$ and for
$
I_-,  I_+,\widetilde{I}_-,  \widetilde{I}_+ \in \Bbb Z,
$
such that
$$
I_+-I_-, \widetilde{I}_+-\widetilde{I}_- \geq 3H,
$$
and for
$$
a \in X_{n}(C)_{(I_-,  I_+ ]}, \quad \widetilde{a} \in X_{n}(C)_{(\widetilde{I}_-, 
 \widetilde{I}_+ ]}, 
$$
such that
$$
a_{(I_-, I_- +  H]} = \tilde{a}_{(\widetilde{I}_-, \widetilde{I}_ ++ H]},\quad
a_{(  I_+  - H,   I_+ ]} = \tilde{a}_{( \widetilde{I}_+ - H,  \widetilde{I}_+]},
$$
one has that $a$ and $\tilde{a}$ have the same context. A subshift $C \subset  \Sigma ^{\Bbb Z}$
has property $(A)$ if there are $H_{n}, n \in  \Bbb N$, such
that $X$ has the properties $(a, n, H_{n}), n \in \Bbb N $.

We  recall  from \cite {Kr} the construction of the associated semigroup. For a  subshift $C \subset \Sigma ^{ \Bbb Z}$ with Property $(A)$ we denote by  $Y(C)$ the set of points in $C$
that are left asymptotic to a  point in $P(X(C))$ and also right asymptotic to a  point 
in $P(X(C))$. Let $y, \widetilde{y}\in Y(C)$, let $y$ be left asymptotic to  $p^{(-)}\in P(X(C))$ and   right asymptotic to $p^{(+)}\in P(X(C))$, let $\widetilde{y}$ be left asymptotic to  $\widetilde{p}^{(-)}\in P(X(C))$ and  right asymptotic to $\widetilde{p}^{(+)}\in P(X(C))$,
 and  let $p^{(-)} \sim_C \widetilde{p}^{(-)}$ and $p^{(+)} \sim_C\widetilde{p}^{(+)}$. The points
$y$ and $\tilde {y}$ are said to be $ \approx_C$-equivalent, if the following holds:
Given that 
$C$ has  properties $(a, n, H_{n}), H_n \geq n,n \in \Bbb N,$ one has
for $n \in \Bbb N$ such that 
$$
{p}^{(-)},\widetilde{p}^{(-)} ,{p}^{(+)},\widetilde{p}^{(+)} 
\in X_n(C),
$$ 
and  for $I_-, I_+, \widetilde
{I}_-, \widetilde {I}_+ \in \Bbb Z,  I_- < I_+, \widetilde {I}_-< \widetilde
{I}_+,$ such that
$$
y_{(- \infty, I_-)} = {p}^{(-)}_{(- \infty, 0)},\quad  y_{(  I_+,\infty)} = {p}^{(+)}_{(  0,\infty)},
$$
$$
\widetilde{y}_{(- \infty, \widetilde{I}_-)} = \widetilde{p}^{(-)}_{(- \infty, 0)},\quad \widetilde{ y}_{(  \widetilde{I}_+\infty)} = \widetilde{p}^{(+)}_{( 
 0,\infty)},
$$
that  
$$ 
a \in C_{[I_- -  3H_n,I_+ +3H_n]},
$$
and
$$
\widetilde {a} \in C_{[\tilde {I}_- -  3H_n,\widetilde {I}_+ + 3H_n]}, 
$$ 
such that
$$
a _{[I_- -H_n ,  I_++H_n]} = y _{[I_- -H_n ,  I_++H_n ]} ,\quad
 \widetilde {a} _{[\widetilde {I}_- -H_n ,  \widetilde {I}_++H_n ]} =
\widetilde{y}_{[\widetilde {I}_- -H_n  ,  \widetilde {I}_++H_n ]},
$$
and
$$
a_{[I_--  3H_n, I_-)} \in X_{n}(C)_{[I_--  3H_n, I_-)} , \quad
\widetilde{a}_{[\widetilde {I}_--  3H_n, \widetilde {I}_-)}
\in X_{n}(C)_{[\widetilde {I}_--  3H_n, \widetilde {I}_-)} ,
$$
$$
 a_{(I_+ , I_++ 3H_n]}  \in X_{n}(C)_{(I_+ , I_++ 3H_n]}  , \quad
  \widetilde{a}_{(\widetilde {I}_+ , \widetilde {I}_++ 3H_n]}  \in 
  X_{n}(C)_{(\widetilde {I}_+ , \widetilde {I}_++ 3H_n]} ,
$$
and
$$
a _{[I_--  3H_n, I_--  2H_n)} = \widetilde {a} _{[\widetilde {I}_- -  3H_n, \widetilde {I}_-  -2 H_{n} )},
$$
$$
a _{(I_++ 2H_{n} ,  I_++ 3H_{n} ]} =
 \widetilde {a} _{(\widetilde {I }_+ +2 H_{n}  , 
\widetilde { I}_+  +3 H_n ]},
$$
have the same context.
The set $[Y(C)]_{\approx}$ carries the
structure of a semigroup: Let $u, v \in Y(C)$, let $u$ be right asymptotic to
$q \in P(X(C))$ and let  $v$ be left asymptotic to
 $r \in P(X(C))$. 
If here 
$q \gtrsim r$, then $[u]_{\approx}[v]_{\approx}$ is  equal to $[y]_{\approx}$,
where $y$ is any point in $Y$, such that there is an
$n \in \Bbb N, $ such that $ q, r \in X_{n}(C),$ together with $I, J  \in \Bbb Z$, and $\widehat {I},  \widehat {J} \in \Bbb Z ,  \widehat {J}- \widehat 
{I}> 3H_n $,  such that
$$
u_{(I, \infty)} = q_{(I, \infty)}, \quad
v_{(-\infty, J]} = r_{(-\infty, J]},
$$  
$$
y_{(- \infty,\widehat{I} + H_{n}]}= u_{(- \infty,I + H_{n}]} ,\quad
y_{ (\widehat{J}  - H_{n}, \infty)}= v_{ (J - H_{n}, \infty)},
$$
and 
$$
y_{(\widehat{I}  , \widehat{J}  ]}\in X_{n}(C)_{(\widehat{I}  , \widehat{J}  ]},
$$
provided that such a point $y$ exists. If such a point $y$ does not exist,  
$[u]_{\approx}[v]_{\approx}$ is  zero.
Also, in the case that
 $q\not \gtrsim r,[u]_{\approx}[v]_{\approx}$ is  zero.
If $\mathcal P_{C}$ is a singleton set $\{ p \}$, then $\mathcal S_{C}$ is a 
monoid (with zero) whose unit is given by $[p]_{\approx_C}$. 

\begin{lemma}
Let C be a topologically transitive subshift, such that
$$
X(C) \neq \emptyset.
$$
Then $Y(C)$ is dense in $ C$.
\end{lemma}
\begin{proof}
It is enough to consider the case that $X_{1}(C) \neq \emptyset.$
 Let $u$ be a forward and backward transitive point of $ C$, and let 
$ 
q \in P(X_1(C))$. Let further $v \in C , K\in \Bbb N$ . Then there are 
$J_{-} , J_{+} , I\in \Bbb Z$, such that 
$$
J_{-} \leq I - K,  \quad J_{+} \geq I + K ,
$$
and such that
$$
u_{I+k} = v_{k}, \qquad -K \leq k \leq K,
$$
and
$$
u_{J_{-}} = u_{J_{+}}  = q_{0}.
$$
One defines a point $y$ in $C$ by
\begin{align*}
&y_{i} = q_{i -J_{-} + I}, \qquad i \leq J_{- }- I,\\
&y_{i} = u_{i + I} , \qquad  J_{-} - I \leq i \leq J_{+ }- I,\\
&y_{i} = q_{i - J_{+} + I},\qquad  J_{+ }- I \leq i,
\end{align*}
and has then
$$
y_{i} = v _{i}    ,\qquad  -K \leq i \leq K,
$$
and
$$
y \in Y(C). \qed
$$
\renewcommand{\qedsymbol}{}
\end{proof}

\begin{lemma}
Let $C \subset \Sigma ^{\Bbb Z}$ be a subshift with 
properties $(a,n, H_n), n \in\Bbb N $. Let $n \in \Bbb N$, let $i_{\circ}, j_{\circ}  \in\Bbb Z, j_{\circ} 
- i_{\circ} > 3H_n$, and let $r \in P(X_n(C)),$ and $ u, v \in 
Y(C),$ 
be 
such that
$$
u_{[i_{\circ}, \infty)} =  r_{[i_{\circ}, \infty)}   ,\quad v_{( - \infty,
j_{\circ}]} =    r_{( - \infty, j_{\circ})}  ,
$$
and 
\begin{align*}
[u]_{\approx(C)}  [v]_{\approx(C)}  \neq 0. \tag {2.1}
\end{align*}
Then 
\begin{align*}
( u _{( - \infty, i_{\circ})},r_{[i_{\circ}, j_{\circ}  )},v_{[j_{\circ}, \infty)}     ) \in C. \tag {2.2}\end{align*}
\end{lemma}
\begin{proof}
By (2.1) there are
$I, J  \in \Bbb Z$, and $\widehat {I},  \widehat {J} \in \Bbb Z ,  \widehat {J}- \widehat 
{I}> 3H_n $,  such that
$$
u_{(I, \infty)} = r_{(I, \infty)}, \quad
v_{(-\infty, J]} = r_{(-\infty, J]},
$$  
$$
y_{(- \infty,\widehat{I} + H_{n}]}= u_{(- \infty,I + H_{n}]} ,\quad
y_{ (\widehat{J}  - H_{n}, \infty)}= v_{ (J - H_{n}, \infty)},
$$
and 
$$
y_{(\widehat{I}  , \widehat{J}  ]}\in X_{n}(C)_{(\widehat{I}  , \widehat{J}  ]}.
$$
Let $k\in \Bbb Z_+,$ be such that $J + k\Pi(r) -I > 3H_n$. By Property $(a, n, H_n)$, then
$$
\Gamma (y_{[\widehat{I}  , \widehat{J})}) =
\Gamma (r_{[{I}  , {J}+ k\Pi(r))}).
$$
For the proof of (2.2) apply again Property $(a, n, H_n)$.
\end{proof}

 \begin{theorem} 
Let $C$ be a topologically transitive subshift with property 
$(A)$.Then the projection of $\mathcal S_{C}$ onto $[\mathcal S_{C}]$ is an isomorphism.
\end{theorem}
\begin{proof}
Let $y, y^\prime \in Y(C)$ be such that 
\begin{align*}
y \not \approx_{C}  y^\prime. \tag {2.3}
\end{align*}
The task is to show that
\begin{align*}
\Gamma([y]_{\approx_C}) \neq \Gamma( [y^\prime]_{\approx_C}). \tag {2.4}
\end{align*}
Let $y$ be left asymptotic to $q \in P(X(C))$ and right asymptotic to $r \in P(X(C))$, and let 
$y^\prime$ be left asymptotic to $q^\prime \in P(X(C))$ and right asymptotic to $r^\prime \in P(X(C)).$ If $q \not \gtrsim_C q^\prime,$ then $[q]_{\approx_C} \not \in 
\Gamma^-([y]_{\approx_C}),$
if $q^\prime \not \gtrsim_C q,$ then $[q^\prime]_{\approx_C} \not \in \Gamma^-([y]_{\approx_C}),$
 if
$r \not \gtrsim_C r^\prime,$ then $[q]_{\approx_C} \not \in \Gamma^+([y]_{\approx_C}),$ and
if $q^\prime \not \gtrsim_C q,$ then $[q^\prime]_{\approx_C}x \not \in 
\Gamma^+([y]_{\approx_C}).$
In the case that $q \sim_C q^\prime$ and $r \sim_C r^\prime$, replacing, if necessary, $y^\prime$ by another representative of 
$[y^\prime]_{\sim_C}$, we can continue with the assumption that $y^\prime$ is left asymptotic to the orbit of $q$ and right asymptotic to the orbit of $r$. Let $n\in \Bbb N$ be such that 
$$
q, r \in P(X_n(C)),
$$
and let $H_n  \in \Bbb N $ be such that the subshift $C$ has Property $(a, n ,H_n)$. Interchanging, if necessary, $y$ and $y^\prime$, one has as
a consequence of Property $(a, n ,H_n)$ and of (2.3), that  there are $K\in \Bbb N$, $K \geq 3H_n\Pi(q)  ,   3H_n\Pi(r),$ and 
$$
(x^-, x^+) \in C_{(-\infty, - K)}\times C_{(K, -\infty)},
$$
such that 
\begin{align*}
y_{(-\infty,-K +3H_n\Pi(q))} =q_{(-\infty, 0)},  \quad
y_{(K-3H_n\Pi(r), \infty)}=r _{[0,\infty)}        ,
\end{align*}
and also
\begin{align*}
y^\prime_{(-\infty,-K +3H_n\Pi(q))} =q_{(-\infty, 0)},  \quad
y^\prime_{(K-3H_n\Pi(r), \infty)}=r _{[0,\infty)} ,   
\end{align*}
and
and such that one has, setting
$$
z = ( x^-   , y_{[-K, K]}  , x^+  ),
\qquad
z^\prime =( x^-   , y^\prime_{[-K, K]}  , x^+  ),
$$
that
$
z \in C,$ but 
\begin{align*}
 z^\prime \not \in C. \tag {2.5}
\end{align*}
By Lemma 2.1 there are points $z^{(m)}\in Y(C), m \in \Bbb N,$ such that
\begin{align*}
z^{(m)}_{[-K-m, K+m]}= z_{[-K-m, K+m]}, \quad m \in \Bbb N. \tag {2.6}
\end{align*}
We construct points $ u^{(m)},v^{(m)}\in Y(C), m \in \Bbb N, $ by setting
$$
 u^{(m)} = ( z^{(m)}_{(-\infty, -K  )} , q_{[-K, \infty)}   ),\qquad 
 v^{(m)} = ( r_{(-\infty, K  ]} , z^{(m)}_{(K, \infty  )} ).
$$
The construction is such that
\begin{align*}
[ u^{(m)}]_{\approx_C} [y]_{\approx_C}  [ v^{(m)}]_{\approx_C} = [z^{(m)}]_{\approx_C}, \qquad m \in \Bbb N. \tag {2.7}
\end{align*}
By (2.5) and (2.6) there is an $m_\circ \in \Bbb N $ such that
\begin{align*}
( z^{(m)}_{(-\infty, -K  )}  ,  y^\prime_{[-K, K]}   ,  z^{(m)}_{(K, \infty  )} )  \not \in C, \qquad m \geq m_\circ. \tag {2.8}
\end{align*}
Therefore
\begin{align*}
[ u^{(m)}]_{\approx_C} [y]_{\approx_C}  [ v^{(m)}]_{\approx_C} =0,\qquad m \geq m_\circ, \tag {2.9}
\end{align*}
since, as a consequence of Lemma 2.2,
$$
[ u^{(m)}]_{\approx_C} [y]_{\approx_C}  [ v^{(m)}]_{\approx_C} \neq 0,
$$
would contradict (2.8). By (2.7) and (2.9), (2.4) is shown.
\end{proof}

 \begin{lemma} 
Let $C$ be a topologically transitive subshift with property 
$(A)$, and let $y \in Y(C)$ be such that
\begin{align*}
[y]_{\approx_C} \in \mathcal M ^{-}_{C}. \tag {2.10}
\end{align*}
Then
$$
y \in X^{-}(C).
$$
\end{lemma}
\begin{proof}
Let $y\in C$ be left asymptotic to $q \in P(X(C))$ and right asymptotic to $r \in P(X(C))$, and let 
$n\in \Bbb N$ be such that 
$$
q, r \in P(X_{n}(C)).
$$
Let $C$  have Property $(a, n ,H_n)$. Let $K \in\Bbb N $ be such that
$$
K > 3H_n\Pi(q),
$$
and
$$
y_{(-\infty, -K)} = q_{(-\infty, 0)}, \qquad y_{[K,\infty)} = r_{[0,\infty)}.
$$
We claim that
$$
y \in X^-_{3K+n}(C).
$$
To justify this claim we let  $ z\in C$
be such that
$$
z_{[-2K, \infty)} = q_{[-K,\infty)},
$$
and we prove that
\begin{align*}
(z_{(-\infty, -K)}   , y_{[-K, \infty)}  )    \in C. \tag  {2.11}
 \end{align*}
By Lemma 2.1 there are points $z^{(m)}\in Y(C), m \in \Bbb N,$ such that
\begin{align*}
z^{(m)}_{[-K-m, \infty)}= z_{[-K-m, \infty)}, \quad m \in \Bbb N. \tag {2.12}
 \end{align*}
 By (2.10) 
 $$
z^{(m)} [y]_{\approx_C}\neq 0,
 $$
and by Lemma 2.2 then
$$
(z^{(m)}_{(-\infty, -K)}   , y_{[-K, \infty)} )     \in C,
$$
and  (2.11) follows now by (2.12).
\end{proof}

\begin{theorem} 
Let $C$ be a topologically transitive subshift with property $(A)$ 
such that
\begin{align*}
\mathcal M ^{-}_{C} \cap  \mathcal M ^{+}_{C}   \neq \emptyset. \tag {2.13}
\end{align*}
Then  $\mathcal M_{C}$  is a monoid (with zero), and
$$
\mathcal M^{-}_{C} \cap  \mathcal M^{+}_{C}   = \{ \bold 1 \}.
$$

\end{theorem}
\begin{proof} 
By (2.13) there is a $y \in Y(C)$  such that 
\begin{align*}
[y]_\approx \in\mathcal M ^{-}_{C} \cap  \mathcal M ^{+}_{C}. \tag {2.14}
\end{align*}
By Lemma 2.4 we have $y \in   X(C).$ Let $y$ be left asymptotic to $q \in P(X(C))$ and right asymptotic to $r \in P(X(C))$. Equation (2.14) implies for all periodic point $p\in X(C)$ that $p \gtrsim_C  q$  and $r  \gtrsim_C  p$. In particular $q\sim_C r$ and it follows that all periodic point $p\in X(C)$ are $\sim_C$-equivalent and $\approx_C$-equivalent to $y$.
\end{proof}

\section{Subshifts constructed from a class of monoids}

 \begin{proposition} 
Let $\mathcal M$ be a monoid, such that
 \begin{align*}
\mathcal  M \setminus \{0 \} =\mathcal  M ^{+}  \mathcal  M ^{-}  , \tag {$\star$}
\end{align*}
and such that
\begin{align*}
\{0, \bold 1\}  \subset( \bigcap_{a^{+}\in \mathcal M^+\setminus \{\bold 1  \}}\mathcal M^- a^{+})
\cap( \bigcap_{a^{-}\in \mathcal M^-\setminus \{\bold 1  \}}a^{-}\mathcal M^+ ) , \tag {$\star$$\star$}
\end{align*}
 Then for $f\in \mathcal M \setminus \{ 0 \}$, the element $a^+\in  \mathcal M^+$  and the  element $a^-\in  \mathcal M^-$ such that
$f = a^+a^-,  $
 are uniquely determined by $f$.
\end{proposition}
\begin{proof}
Let for $f\in \mathcal M \setminus \{ 0 \}$ be given $a^+,b^+ \in  \mathcal M^+ $ and
$a^-,b^- \in  \mathcal M^- $ such that
\begin{align*}
f =a^+a^- =b^+b^-.  \tag {3.1}
\end{align*}
By ($\star$$\star$) there are $ c^+,d^+\in  \mathcal M^+,  $ such that
\begin{align*}
a^-c^+ = b^-d^+=\bold 1. \tag {3.2}
\end{align*}
Then
\begin{align*}
b^-c^+\in  \mathcal M^+.  \tag {3.3}
\end{align*}
Otherwise, there would be by ($\star$) a $g^+ \in \mathcal M^+$ and a
$g^- \in \mathcal M^-\setminus \{  \bold 1\},$
such that
\begin{align*}
b^-c^+ =g^+g^-, \tag {3.4}
\end{align*}
and by ($\star$$\star$) there would then be an $h^+ \in \mathcal M^+ \setminus \{  \bold 1\}$ such that
\begin{align*}
g^-h^+= 0. \tag {3.5}
\end{align*}
One derives now a contradiction from (3.1): On the one hand one has by (3.2)  that
$$
a^+a^-c^+ h^+ = a^+ h^+ \in \mathcal M^+,
$$
and on the other hand one has by (3.4) and (3.5)
that
$$
b^+b^-c^+ h^+ =b^+g^+g^- h^+ = 0.
$$
From (3.1)  and (3.2) 
\begin{align*}
a^+ = b^+ b^- c^+. \tag {3.6}
\end{align*}
and by the same argument
\begin{align*}
b^+ = a^+ a^- d^+. \tag {3.7}
\end{align*}
From (3.6) and (3.7)
$$
b^+ =b^+ b^- c^+a^- d^+,
$$
which implies as a consequence of ($\star$$\star$) that
$$
b^- c^+a^- d^+ = \bold 1,
$$
and in the same way it follows from (3.3)  that
$$
a^- d^+ = \bold 1.
$$
From  (3.7) one has then that  $b^+ = a^+$ and from ($\star$$\star$) also that $b^- = a^-$.
\end{proof}

\begin{theorem} 
Let $\mathcal M$ be a monoid,  
such that ($\star$) and ($\star$$\star$) hold.
and let $\Delta$ be a finite generating set of $\mathcal M$. Then the subshift $C(\Delta)$ has Property $(A)$, and
$$
\mathcal S_{C(\Delta) } = [\mathcal M].
$$
\end{theorem}
\begin{proof}
For $f \in \mathcal M \setminus \{0\}$ choose $I_f\in \Bbb N,$ and $\delta _i^{(f)}\in \Delta, 1 \leq i \leq I_f,$ such that
$$
\prod_{1 \leq i \leq I_f }\delta^{(f)} _i= f,
$$
and set
$$
w^{(f)} = (\delta^{(f)} _i  )_{1 \leq i \leq I_f }.
$$

Let $H, I \in \Bbb N, I \geq 2H$, and let $u\in C(\Delta)_{[1, I]}$ be such that
\begin{align*}
u_{(H, I]} \in \omega^+(u_{(1, H]}   ), \qquad u_{[1, I-H]} \in \omega^-(u_{(I-H, I]}   ). \tag {3.8}
\end{align*}
We apply Proposition (3.1).
We let
$$
a^-(-), c^-,a^-(+)  \in \mathcal M^-, \qquad a^+(-), c^+,a^+(+)  \in \mathcal M^+,
$$
be 
given by
$$
\prod_{1\leq i \leq H}u_i = a^+(-)a^-(-)  ,
\prod_{H < i \leq I - H}u_i = c^+c^- ,
\prod_{I - H < i \leq I}u_i=  a^+(+)a^-(+) ,
$$
and let then 
$b^- \in \mathcal M^-, b^+ \in \mathcal M^+$ 
given by
$$
a^-(-) c^+c^-  a^+(+)  =  b^+b^-.
$$
It is $b^+ = \bold 1$, for otherwise there  would exist by 
($\star$$\star$), a $ g^-\in \mathcal M^-$ and an $h^- \in \mathcal M^-$ such that $g^-a^+(-)  = 1,$ and $h^-b^+  = 0,$
and then 
$$
(w^{(h^-)}   , w^{(g^-)}   ) \in \Gamma^-(u_{[1,H]}), 
$$
but
$$
(w^{(h^-)}   , w^{(g^-)}   ) \notin \Gamma^-(u),
$$
contradicting (3.8). By the symmetric argument one has that $b^- = \bold 1$. 
It follows for $ v(-)\in C(\Delta)_{(- I_-, 1]},  v(+)  \in C(\Delta)_{(I, I+I_+]},  I_-,I_+ \in \Bbb N,$ that $(v(-),v(+)) \in \Gamma(u)$ if and only if
$$
\left(\prod_{- I_-< i \leq 1} v_i(-)\right) a^+(-)a^-(+)
\left(\prod_{ I < i \leq I+I_+}v_i(+)\right)\neq 0,
$$
and this shows that 
$C(\Delta)$ has Property $(a, H, H)$. It also follows that a periodic point $p$ is in $X(C(\Delta))$ if and only of $p$ has a period $\Pi$ such that with 
$a^- \in \mathcal M^-,a^+ \in \mathcal M^+,$ given by
$$
\prod_{1\leq i \leq \Pi}p_i = a^+a^-  ,
$$
one has  $a^-a^+=  \bold 1$.

We show that all periodic points in $X(C(\Delta))$ are $\sim_{C}$-equivalent. For this let $q$ and $r$ be periodic points in $X(C(\Delta))$, and let $\Pi_q$ be a period of $q$, and  $\Pi_rq$  a period of $r$, such that with $a^-(q), a^-(r) \in  \mathcal M^-$ and $a^+(q), a^+(r) \in  \mathcal M^+,$ 
given by
$$
\prod_{1\leq i \leq \Pi_q}q_i = a^+(q)a^-(q),  \qquad \prod_{1\leq i \leq \Pi_r}r_i = a^+(r)a^-(r) ,
$$
one has
$$ 
 a^-(q)a^+(q)  = \bold 1, \qquad   a^-(r)a^+(r)  = \bold 1.
$$
By ($\star$$\star$) there are $b^-\in  \mathcal M^-$ and $b^+\in  \mathcal M^+$ such that 
$$
a^-(q)b^+ = \bold 1,  \qquad b^-a^+(r) = \bold 1.
$$
The point $x \in C(\Delta)$, that is given by setting
$$
x_{[1, I_{b^+b^-}]} = w^{(b^+b^-)},
$$
and
$$
x_i = \begin{cases} q_i, &\text{if $i < 0,$} \\
r_{i-I_{b^-b^+}} , &\text{if $i > I_{b^-b^+},$}
\end{cases}
$$
is in $X(C)$.

For $f\in \mathcal M$ let $x^{(f)}$ be the point in $Y(C(\Delta))$ that is given by setting
\begin{align*}
&x^{(f)}_{(-(k+1)I_{\bold 1}, kI_{\bold 1}]} = w^{(\bold 1)}, \qquad\qquad \thinspace k \in \Bbb N,
\\
&x^{(f)}_{[1,I_{f}]}=w^{(f)},
\\
&x^{(f)}_{(I_f+kI_{\bold 1}, I_f+(k+1)I_{\bold 1})} = w^{(\bold 1)}, \qquad k \in \Bbb Z_+.
\end{align*}
The $\approx_C$-class of $x^{(f)}$ does not depend on the choice of $w^{(f)}$, and  a homomorphism $\varphi$ of $\mathcal S_{C(\Gamma)}$  onto $[\mathcal M]$ is obtained by setting
$$
\varphi(f) =[x^{(f)}]_{\approx_C}, \qquad f\in \mathcal M.
$$ 
The homomorphism $\varphi$ is injective since for $f, f^\prime\in \mathcal M,$
$$
\Gamma(f) =\Gamma(f^\prime) ,
$$
if and only if
$$
\Gamma( x^{(f)}_{ [-kI_{\bold 1}  , I_{f}+kI_{\bold 1}   ]}   ) =
\Gamma( x^{(f^\prime)}_{ [-kI_{\bold 1}  ,  I_{f^\prime}+kI_{\bold 1}   ]}), \quad
k \in \Bbb Z_+. \qed
$$
\renewcommand{\qedsymbol}{}
\end{proof}

We note at this point that the subshifts $C(\Delta)$ are coded systems in the sense of \cite{BH}.
Also compare to the definition of the subshifts $C(\Delta)$ how sofic systems were originally defined \cite {W}.

 \begin{proposition} 
 Let $\mathcal M$ be a monoid, 
 such that ($\star$) and ($\star$$\star$) hold.
Let the mappings
$$
 a^{+} \to \Gamma^{-}(a^{+})\ ( a^{+} \in \mathcal M^{+}), \quad
 a^{-} \to \Gamma^{+}(a^{-})\ ( a^{-} \in \mathcal M^{-}),
$$
be injective. Then the projection of $\mathcal M$ onto $[\mathcal M]$ is an
isomorphism.
\end{proposition}
\begin{proof}
 Let $a ,\widetilde{a} \in \mathcal M , a \neq \widetilde{a}$, and let 
 $b^{+}, \widetilde{b}^{+}\in  \mathcal M^+, b^{-}, \widetilde{b}^{-}\in  \mathcal M^-$
be  given by
$$
 a = b^{+} b^{-}   , \widetilde{a} = \widetilde{b}^{+} \widetilde{b}^{-}    , $$
and consider the case that
$$
b^{+}   \neq \widetilde{b}^{+}.
$$
Let $c^{-} \in \mathcal M^{-}$ be such that
$$
c^{-}  \in \Gamma^{-}(b^{+} ), c^{-}  \notin
\Gamma^{-}( \widetilde{b}^{+}) .
$$
Then
$$
 (c^{-} ,\bold 1) \in  \Gamma(a),\quad (c^{-} ,\bold 1)    \notin 
\Gamma( \widetilde{a} ) .\qed
$$  
\renewcommand{\qedsymbol}{}
\end{proof}

\section{A class of examples}

Consider the case of a  monoid $\mathcal M$ that satisfies ($\star$) and ($\star$$\star$) where the submonoids $\mathcal M^-$ and $\mathcal M^+$ are finitely and freely generated. With
$N_- ,N_+  \in \Bbb N$ and a generating set
$\Delta_- =\{   \delta^-(n_-): 0 \leq n_- \leq N_-     \}$ of $\mathcal M^-$
and a 
generating set
$\Delta_ +=\{   \delta^+(n_+): 0 \leq n_+ \leq N_+    \}$ of $\mathcal M^+$
there are then relations
\begin{multline*}
\delta^-(n_-) \delta^+(n_+) =\tag {4.1}\\
\left ( \prod_{I_+( n_- ,  n_+) \geq i_{+, n_- ,  n_+}>0}
    \delta^+(i_{+, n_- ,  n_+}   )     \right) 
    \left( \prod_{0< i_{+, n_- ,  n_+}\leq I_+( n_- ,  n_+))}
    \delta^-( i_{+, n_- ,  n_+}  )     \right ),  \\
    0 \leq n_- \leq N_-,0 \leq n_+ \leq N_+.
\end{multline*}
This suggests  to consider the subclass of the class of monoids (with zero) that satisfy ($\star$) and  ($\star$$\star$) which contains the monoids  (with zero) that are generated by $\Delta_-\cup\Delta_+   $ with the  relations that are expressed by (4.1), in which case (4.1) is necessarily a confluent string-rewriting system. Conversely, given a confluent string-rewriting system  (4.1) there are procedures to decide if it defines a monoid $\mathcal M$ (with zero) that satisfies  ($\star$$\star$), and if the projection of $\mathcal M$  onto  $[\mathcal M ]$ is an isomorphism \cite[Section 4.3]{BO}.

For the examples, which we are going to write as string-rewriting systems, it can be  verified by inspection that they describe monoids (with zero) $\mathcal M$ such that  ($\star$$\star$) holds, and by means of Proposition 3.3 it can be  verified that the projection of $\mathcal M$  onto  $[\mathcal M ]$ is an isomorphism.

A prototype example is the Dyck inverse monoid (the polycyclic monoid \cite {NP}), that is given by
\begin{align*}
& \Delta_- =  \{ \lambda, \lambda^{\prime} \}   ,\Delta_+= \{ \rho,  \rho^{\prime}   
\}, \\& \lambda \rho ,\lambda^{\prime}  \rho^{\prime}   \to  \bold 
1,     \\
&\lambda\rho^{\prime},
 \lambda^{\prime}\rho \to 0.
\end{align*}
Examples where $\lambda\rho \in \{ 0, \bold 1 \}, \lambda \in 
\Delta_-,
\rho \in \Delta_+ $ (compare here [Ke]) are
\begin{align*}
& \Delta_- = \{ \lambda, \lambda^{\prime} ,\lambda^{\prime\prime} \}  , 
       \Delta_+= \{\rho, 
\rho^{\prime}  , \rho^{\prime\prime} \},\\
&  \lambda \rho^{\prime} ,\lambda \rho,\lambda^{\prime}  
\rho^{\prime} ,
\lambda^{\prime\prime}  \rho^{\prime\prime}   \to  \bold 1,\\
&\lambda ^{\prime}\rho,\lambda \rho^{\prime\prime},\lambda^{\prime}
\rho^{\prime\prime},\lambda^{\prime\prime} \rho,\lambda
^{\prime\prime}\rho^{\prime} \to 0,
\end{align*}
and
\begin{align*}
& \Delta_- = \{ \lambda, \lambda^{\prime} ,\lambda^{\prime\prime}\}  ,  \Delta_+ = 
\{\rho, \rho^{\prime}  , \rho^{\prime\prime} \}, \\
& \lambda \rho ,\lambda^{\prime}  \rho,\lambda^{\prime}  
\rho^{\prime} ,\lambda^{\prime\prime}  \rho^{\prime\prime}  \to  
\bold 1,\\
& \lambda  \rho^{\prime}, \lambda \rho^{\prime\prime}, 
\lambda^{\prime}\rho^{\prime\prime} ,\lambda^{\prime\prime} \rho,\lambda^{\prime\prime} \rho^{\prime} \to  0.
\end{align*}
Another example is
\begin{align*}
&  \Delta_- =  \{ \lambda, \lambda^{\prime} ,\lambda^{\prime\prime}  \} , 
 \Delta_+ =\{\rho, 
\rho^{\prime}  , \rho^{\prime\prime} \}, \\&  \lambda \rho 
,\lambda^{\prime} 
\rho^{\prime} , \lambda^{\prime\prime}  \rho^{\prime\prime}   \to  
\bold 1,\\
  &  \lambda    \rho^{\prime\prime} ,\lambda^{\prime}\rho, 
\lambda^{\prime\prime}\rho^{\prime} \to 0.\\
&\lambda \rho^{\prime} \to   \lambda,   \\
& \lambda^{\prime}\rho^{\prime\prime} \to  \lambda^{\prime} ,  \\
&\lambda^{\prime\prime}\rho^{\prime} \to  \lambda^{\prime\prime}. 
\end{align*}

\medskip 

Toshihiro Hamachi

Faculty of Mathematics

Kyushu University

744 Motooka, Nishi-ku

Fukuoka 819-0395, 

Japan

hamachi@math.kyushu-u.ac.jp

\bigskip

Wolgang Krieger

Institute for Applied Mathematics

University of Heidelberg

Im Neuenheimer Feld 294

69120 Heidelberg

Germany

krieger@math.uni-heidelberg.de

 \end{document}